\documentclass[preprint,sort&compress,3p]{elsarticle}
\usepackage[utf8]{inputenc}
\usepackage{bbm}
\usepackage{amsmath}
\usepackage{amssymb}
\usepackage{amsthm}
\usepackage{graphicx}
\usepackage{color}
\usepackage{enumerate}

\def\dualsymb{\circ}
\def\transposesymb{t}
\def\conjsymb{\ast}

\newcommand{\setC}{\mathbbm{C}}

\newcommand{\bounded}[1]{\mathcal{B}\left(#1\right)}

\newcommand{\boundedconj}[1]{\bounded{#1}^{\conjsymb}}
\newcommand{\bbounded}[2]{\mathcal{B}\left(#1,#2\right)}
\newcommand{\hilbertspaceone}{\mathcal{K}}
\newcommand{\hilbertspacetwo}{\mathcal{H}}

\newcommand{\bk}{\bounded{\hilbertspaceone}}
\newcommand{\bh}{\bounded{\hilbertspacetwo}}

\newcommand{\mappingcone}{\mathcal{C}}
\newcommand{\innerpr}[2]{\left<#1,#2\right>}

\newcommand{\Tr}{\mathop{\textnormal{Tr}}}
\newcommand{\Ad}{\mathop{\textnormal{Ad}}}
\newcommand{\rk}{\mathop{\textnormal{rk}}}
\newcommand{\supp}{\mathop{\textnormal{supp}}}

\newcommand{\conj}[1]{#1^{\conjsymb}}
\newcommand{\dual}[1]{#1^{\dualsymb}}
\newcommand{\ddual}[1]{#1^{\dualsymb\dualsymb}}
\newcommand{\convhull}{\mathop{\textnormal{convhull}}}
\newcommand{\Id}{\mathop{\textnormal{id}}}
\newcommand{\One}{\mathbbm{1}}
\newcommand{\proj}[1]{p_{#1}}
\newcommand{\diad}[1]{\left|#1\right>\left<#1\right|}
\newcommand{\ddiad}[2]{\left|#1\right>\left<#2\right|}
\newcommand{\Choimatr}[1]{C_{#1}}
\newcommand{\Choimatrplus}[1]{\Choimatr{#1}^+}
\newcommand{\Choimatrminus}[1]{\Choimatr{#1}^-}
\newcommand{\matrices}[2]{M_{#1}\left(#2\right)}
\newcommand{\transpose}{\transposesymb}

\newcommand{\Pmaps}{\mathcal{P}}
\newcommand{\kPmaps}[1]{\Pmaps_{#1}}
\newcommand{\SPmaps}{\mathcal{SP}}
\newcommand{\kSPmaps}[1]{\SPmaps_{#1}}
\newcommand{\CPmaps}{\mathcal{CP}}

\newcommand{\Pmapsb}[1]{\Pmaps\left(#1\right)}
\newcommand{\kPmapsb}[2]{\kPmaps{#1}\left(#2\right)}
\newcommand{\SPmapsb}[1]{\SPmaps\left(#1\right)}
\newcommand{\kSPmapsb}[2]{\kSPmaps{#1}\left(#2\right)}
\newcommand{\CPmapsb}[1]{\CPmaps\left(#1\right)}

\newcommand{\Bplusbb}[2]{\mathcal{B}\left(#1,#2\right)^+}
\newcommand{\Kpositive}[1]{\mathcal{P}_{#1}}

\newcommand{\sk}[1]{{\mathcal S}_{#1}}
\newcommand{\norm}[1]{\left\|#1\right\|} 
\newcommand{\abs}[1]{\left|#1\right|}
\newcommand{\kfnorm}[2]{\left\|#2\right\|_{\left(#1\right)}}
\newcommand{\kkfnorm}[1]{\kfnorm{k}{#1}}
\newcommand{\hsnorm}[1]{\norm{#1}_{\textnormal{HS}}}
\newcommand{\hsprod}[2]{\innerpr{#1}{#2}_{\textnormal{HS}}}

\def\Cdot{\,\cdot\,}

\newtheorem{theorem}{Theorem}
\newtheorem{lemma}{Lemma}

\newtheorem{definition}{Definition}
\newtheorem{proposition}{Proposition}

\newtheorem{corollary}{Corollary}

\begin{document}

\title{Choi matrices, norms and entanglement associated with positive maps on matrix algebras}

\author{{\L}ukasz~Skowronek}
\ead{lukasz.skowronek@uj.edu.pl}
\address{Instytut Fizyki im. Smoluchowskiego, Uniwersytet
Jagiello{\'n}ski, Reymonta 4, 30-059 Krak{\'o}w, Poland 
}

\author{Erling St\o rmer}
\ead{erlings@math.uio.no}
\address{Department of Mathematics, 
	University of Oslo, P.O. Box 1053 Blindern,
	NO-0316 Oslo, Norway}

\begin{abstract}
We study positive maps of $\bk$ into $\bh$ for finite-dimensional Hilbert spaces $\hilbertspaceone$ and $\hilbertspacetwo$. Our main emphasis is on how Choi matrices and estimates of their norms with respect to mapping cones reflect various properties of the maps. Special attention will be given to entanglement properties and $k$-positive maps, in particular tensor products of $2$-positive maps. The latter problem is directly related to the question of $n$-copy distillability of quantum states, for which we obtain a partial result.
\end{abstract}


\maketitle

\section{Introduction}
Positive maps between matrix algebras have been studied in quantum information theory for at least fifteen years  and in operator algebras since the $1950$'s. In the study several positivity conditions have been introduced, especially completely positive maps and $k$-positive maps have been of importance. Positive maps are completely characterized by their Choi matrices, which are also recognized as entanglement witnesses for states. The Choi matrix for a positive map is self-adjoint, hence is the orthogonal difference of two positive matrices, which implies that each positive map is a difference of two completely positive maps, its negative and its positive part. In the paper we shall see that the negative part contains important information.

A central class of maps are those of the form $\phi_{\lambda}\left(x\right)=\Tr\left(x\right)\One-\lambda\psi\left(x\right)$, where $\psi$ is a completely positive map and $\Tr$ is the usual trace. We shall give a characterization of such a map to be so-called $\mappingcone$-positive, where $\mappingcone$ stands for an arbitrary mapping cone. In the particular case when $\mappingcone$ is equal to the set of $k$-positive maps and $\psi=\Ad_V$, we recover a recent result of Chru{\'s}ci{\'n}ski and Kossakowski \cite{ref.ChK09}, thus putting it in a more general framework. Furthermore, we use maps like $\phi_{\lambda}$ to relate $\mappingcone$-positivity of maps to the properties of the support of the negative part of their Choi matrices. Using our results we also provide conditions for the tensor product $\phi_{\lambda}\otimes\phi_{\lambda}$ to be $2$-positive. This latter result is related to a work by Di Vicenzo et al. \cite{ref.DiVicenzo00}, in which they show that it is sufficient to study maps of the form $\phi_{\lambda}\left(x\right)=\Tr\left(x\right)-\lambda\Ad_V$ and their tensor products to answer the fundamental question about the existence of NPT bound entanglement.

\section{Norms and positivity conditions}

We denote by $\bh$ the algebra of bounded operators on a Hilbert space $\hilbertspacetwo$. Throughout the paper all Hilbert spaces will be finite-dimensional, so of the form $\setC^d$, $d$ the dimension of $\hilbertspacetwo$. Let $\hilbertspaceone$ be another Hilbert space. Then we denote by $\bbounded{\bk}{\bh}$ (resp. $\Bplusbb{\bk}{\bh}$) the linear (resp. positive linear) maps of $\bk$ into $\bh$. Here positive means that positive operators are mapped to positive operators. A map $\phi$ is called {\it $k$-positive} if $\phi\otimes\Id_{\matrices{k}{\setC}}$, considered as a map in $\bbounded{\bounded{\hilbertspaceone\otimes\setC^k}}{\bounded{\hilbertspacetwo\otimes\setC^k}}$ is positive, and {\it completely positive} if it is $k$-positive for all $k\geqslant 1$. In the latter case $\phi$ is a sum of maps of the form $\Ad_V$, which sends $x\mapsto Vx\conj{V}$, where $V:\hilbertspaceone\rightarrow\hilbertspacetwo$, with $V$ linear \cite{ref.Choi75, ref.Kraus}. If $\phi=\sum_i\Ad_{V_i}$, where the rank $\rk V_i\leqslant k$ for all $i$, then $\phi$ is called {\it $k$-superpositive}, in particular if $k=1$, then $\phi$ is superpositive \cite{ref.Ando04}, or entanglement breaking \cite{ref.HSR03}.

Let $\Pmapsb{\hilbertspacetwo}:=\Bplusbb{\bh}{\bh}$. A closed cone $\mappingcone$ in $\Pmapsb{\hilbertspacetwo}$ is called a {\it mapping cone} \cite{ref.St86} if $\alpha\circ\phi\circ\beta\in\mappingcone$ whenever $\alpha,\beta\in\CPmapsb{\hilbertspacetwo}$ - the completely positive maps in $\Pmapsb{\hilbertspacetwo}$. Equivalently, by the above decomposition of completely positive maps, $\phi\in\mappingcone$ iff $\Ad_U\circ\phi\circ\Ad_V\in\mappingcone$ for all $U,V\in\bh$. We shall denote the set of $k$-positive maps in $\Pmapsb{\hilbertspacetwo}$ by $\kPmapsb{k}{\hilbertspacetwo}$, or  
just $\kPmaps{k}$ when the context is clear. Superpositive maps in $\Pmapsb{\hilbertspacetwo}$ will be denoted with $\SPmapsb{\hilbertspacetwo}$ and $k$-superpositive maps by $\kSPmapsb{k}{\hilbertspaceone}$.  It is a simple exercise to check that the cones $\Pmapsb{\hilbertspacetwo}$, $\kPmapsb{k}{\hilbertspacetwo}$, $\CPmapsb{\hilbertspacetwo}$, $\kSPmapsb{k}{\hilbertspacetwo}$, $\SPmapsb{\hilbertspacetwo}$ of positive, $k$-positive, completely positive, $k$-superpositive and superpositive maps (resp.), are examples of mapping cones. They are also {\it symmetric mapping cones}, i.e. $\phi\in\mappingcone$ implies that the maps $\transpose\circ\phi\circ\transpose$ with $t$ the transpose operation and $\conj{\phi}$, defined by $\Tr\left(\phi\left(a\right)b\right)=\Tr\left(a\conj{\phi}\left(b\right)\right)$, belong to $\mappingcone$. If $\mappingcone$ is a mapping cone in $\Pmapsb{\hilbertspacetwo}$ then one of us (E.S., \cite{ref.St86}) introduced a positivity property for maps in $\Bplusbb{\bk}{\bh}$ called $\mappingcone$-positivity. It was shown in \cite{ref.St09mappingcones} that if $\mappingcone$ is symmetric, then the cone of $\mappingcone$-positive maps, denoted by $\Kpositive{\mappingcone}$ in the sequel, are the maps in the closed cone in $\Bplusbb{\bk}{\bh}$ generated by maps of the form $\alpha\circ\psi$, where $\alpha\in\mappingcone$, and $\psi$ is a completely positive map of $\bk$ into $\bh$. In particular, if $\hilbertspacetwo=\hilbertspaceone$, $\Kpositive{\mappingcone}=\mappingcone$. It follows that $\phi\in\Pmapsb{\hilbertspacetwo}$ belongs to $\Kpositive{\kPmaps{k}}$ iff $\phi\in\kPmaps{k}$, and $\phi\in\Kpositive{\kSPmaps{k}}$ iff $\phi\in\kSPmaps{k}$.

Let $\left(e_{ij}\right)$ be a complete set of matrix units for $\bk$. Then there is a one-to-one correspondence between maps $\phi\in\bbounded{\bk}{\bh}$ and operators in $\bounded{\hilbertspaceone\otimes\hilbertspacetwo}$ ($=\bk\otimes\bh$), given by $\phi\mapsto\Choimatr{\phi}=\sum_{i,j}e_{ij}\otimes\phi\left(e_{ij}\right)$. The operator $\Choimatr{\phi}$ is called  the {\it Choi matrix} for $\phi$ \cite{ref.Choi75}, and the map $J:\phi\mapsto\Choimatr{\phi}$ is sometimes called the {\it Jamiołkowski-Choi isomorphism} \cite{ref.J72}. It was shown by Choi that $\phi$ is completely positive iff $\Choimatr{\phi}$ is a positive matrix. More generally, $\phi$ is $k$-positive iff $\Choimatr{\phi}$ is {\it $k$-block positive} (cf. \cite{ref.SSZ09}). In particular, $\phi$ is positive iff $\Tr\left(\Choimatr{\phi}a\otimes b\right)\geqslant 0$ for all positive $a\in\bk$, $b\in\bh$.

Given a cone $\mappingcone$ in $\Bplusbb{\bk}{\bh}$, one defines its dual cone $\dual{\mappingcone}$ by the formula
\begin{equation}\label{dualdef}
 \dual{\mappingcone}=\left\{\psi\in\bbounded{\bk}{\bh}\vline\Tr\left(\Choimatr{\phi}\Choimatr{\psi}\right)\geqslant 0\,\forall_{\phi\in\mappingcone}\right\}.
\end{equation}

If $\mappingcone$ is a mapping cone, we denote 
\begin{equation}\label{defsk}
 \sk{\mappingcone}:=\left\{\rho\in\boundedconj{\hilbertspaceone\otimes\hilbertspacetwo}\,\vline\,\rho=\Tr\left(\Choimatr{\psi}\Cdot\right),\Tr\left(\Choimatr{\psi}\right)=1,\psi\in\dual{\Kpositive{\mappingcone}}\right\}.
\end{equation}
For example, $\sk{\Pmaps}$ denotes the set of separable states on $\bounded{\hilbertspaceone\otimes\hilbertspacetwo}$. We can now define a norm on $\bounded{\hilbertspaceone\otimes\hilbertspacetwo}$ by
\begin{equation}\label{defnormsk}
 \norm{A}_{\sk{\mappingcone}}:=\sup_{\rho\in\sk{\mappingcone}}\abs{\rho\left(A\right)}.
\end{equation}
There is a corresponding norm on $\bbounded{\bk}{\bh}$, defined by
\begin{equation}\label{defnormk}
 \norm{\psi}_{\mappingcone}=\sup_{\phi\in J^{-1}\left(\sk{\mappingcone}\right)}\abs{\Tr\left(\Choimatr{\phi}\Choimatr{\psi}\right)}=\norm{\Choimatr{\psi}}_{\sk{\mappingcone}},
\end{equation}
where we identify a linear functional with its density operator.
The norm properties $\norm{\lambda\psi}_{\mappingcone}=\abs{\lambda}\norm{\psi}_{\mappingcone}$ and $\norm{\phi+\psi}_{\mappingcone}\leqslant\norm{\phi}_{\mappingcone}+\norm{\psi}_{\mappingcone}$ are immediate from definition. To show that $\norm{\psi}_{\mappingcone}=0$ implies $\psi=0$, note that by \cite[Lemma 2.4]{ref.St86} each mapping cone contains the superpositive maps.  Since the composition of a superpositive map and 
a positive map is completely positive, the superpositive maps of  $\bk$ into $\bh$ belong to the dual 
cone $\dual{\Kpositive{\mappingcone}}$ of $\Kpositive{\mappingcone}$ by \cite[Thm. 1]{ref.St09dual}.  Thus $\sk{\mappingcone}$ contains all states with 
density operators  corresponding to maps in $\SPmaps$.  Since these states 
form a separating family of states, $\Choimatr{\psi}=0$. Hence $\psi = 0$.

Note that for $\mappingcone=\kPmaps{k}$, $k=1,\ldots,d$, $d=\dim\hilbertspacetwo$ and $A$ normal, the norm $\norm{A}_{\sk{\mappingcone}}$ reduces to the Schmidt norm $\norm{A}_{S\left(k\right)}$ introduced by Johnston and Kribs \cite{ref.JK09}.

If $\phi$ is a positive map of $\bk$ into $\bh$, then $\Choimatr{\phi}$ is a self-adjoint operator in $\bounded{\hilbertspaceone\otimes\hilbertspacetwo}$, hence is a difference $\Choimatr{\phi}=\Choimatrplus{\phi}-\Choimatrminus{\phi}$ of two positive operators $\Choimatrplus{\phi}$ and $\Choimatrminus{\phi}$ such that $\Choimatrplus{\phi}\Choimatrminus{\phi}=0$. Let $\phi^+=J^{-1}\left(\Choimatrplus{\phi}\right)$, $\phi^-=J^{-1}\left(\Choimatrminus{\phi}\right)$. Since $\Choimatrplus{\phi}$ and $\Choimatrminus{\phi}$ are positive, $\phi^+$ and $\phi^-$ are completely positive by the Choi theorem \cite{ref.Choi75}.

\begin{proposition}\label{normCplusCminus}
 Let $\phi$ belong to a mapping cone $\mappingcone\supset\CPmapsb{\hilbertspacetwo}$. With the above notation
\begin{equation}\label{normineqphi}
 \norm{\phi_+}_{\mappingcone}\geqslant\norm{\phi_-}_{\mappingcone}
\end{equation}
or equivalently, $\norm{\Choimatrplus{\phi}}_{\sk{\mappingcone}}\geqslant\norm{\Choimatrminus{\phi}}_{\sk{\mappingcone}}$.
\begin{proof}
 If $\psi\in J^{-1}\left(\sk{\mappingcone}\right)\subset\dual{\Kpositive{\mappingcone}}=\mappingcone^{\circ}$, we have
\begin{equation}\label{prop1proof1}
 0\leqslant\Tr\left(\Choimatr{\phi}\Choimatr{\psi}\right)=\Tr\left(\Choimatrplus{\phi}\Choimatr{\psi}\right)-\Tr\left(\Choimatrminus{\phi}\Choimatr{\psi}\right).
\end{equation}
Thus
\begin{equation}\label{prop1proof2}
 \norm{\phi^+}_{\mappingcone}-\sup_{\psi\in J^{-1}\left(\sk{\mappingcone}\right)}\Tr\left(\Choimatrminus{\phi}\Choimatr{\psi}\right)\geqslant 0.
\end{equation}
Since $\mappingcone\supset\CPmaps$, $\dual{\mappingcone}\subset\dual{\CPmaps}=\CPmaps$. Thus $\Tr\left(\Choimatrminus{\phi}\Choimatr{\psi}\right)\geqslant 0$ for all $\psi\in J^{-1}\left(\sk{\mappingcone}\right)$. Therefore
\begin{equation}\label{prop1proof3}
 \norm{\phi^+}_{\mappingcone}\geqslant\sup_{\psi\in J^{-1}\left(\sk{\mappingcone}\right)}\abs{\Tr\left(\Choimatrminus{\phi}\Choimatr{\psi}\right)}=\norm{\phi^-}_{\mappingcone}.
\end{equation}
Since $\norm{\Choimatrplus{\phi}}_{\sk{\mappingcone}}=\norm{\phi^+}_{\mappingcone}$, and the same for $\phi^-$ and $\Choimatrminus{\phi}$ the proof is complete.
\end{proof}
\end{proposition}
Let us now consider maps $\phi_{\lambda}$ of the form $\phi_{\lambda}\left(a\right)=\Tr\left(a\right)\One-\lambda\phi\left(a\right)$, $a\in\bk$, where $\phi$ is a completely positive map of $\bk$ into $\bh$, and $\lambda\in\left[0,+\infty\right)$. For notational convenience we identify a linear functional $\omega$ on $\bk$ with the map $a\mapsto\omega\left(a\right)\One$, where $\One$ is the identity in $\bh$. Thus we shall write $\phi_{\lambda}=\Tr-\lambda\phi$.
\begin{proposition}\label{normphilambda}Let $\mappingcone$ be a symmetric mapping cone on $\hilbertspacetwo$ and $\phi_{\lambda}=\Tr-\lambda\phi$ as above. Then $\phi_{\lambda}$ is $\mappingcone$-positive iff
\begin{equation}\label{ineqphilambda}
 \norm{\phi}_{\mappingcone}\leqslant\frac{1}{\lambda}.
\end{equation}
\begin{proof}
Since $\ddual{\Kpositive{\mappingcone}}=\Kpositive{\mappingcone}$ (see e.g. \cite[Thm. 6]{ref.St09dual}), $\phi_{\lambda}$ is $\mappingcone$-positive iff $\Tr\left(\Choimatr{\phi_{\lambda}}\Choimatr{\psi}\right)\geqslant 0\,\forall_{\psi\in\dual{\Kpositive{\mappingcone}}}$, iff $\Tr\left(\Choimatr{\phi_{\lambda}}\Choimatr{\psi}\right)\geqslant 0\,\forall_{\psi\in J^{-1}\left(\sk{\mappingcone}\right)}$, hence iff $\inf_{\psi\in J^{-1}\left(\sk{\mappingcone}\right)}\Tr\left(\Choimatr{\phi_{\lambda}}\Choimatr{\psi}\right)\geqslant 0$. Now $\Choimatr{\Tr}=\sum_{i,j}e_{ij}\otimes\Tr\left(e_{ij}\right)=\sum_{i}e_{ii}\otimes\One=\One$. Thus $\phi$ is $\mappingcone$-positive iff
\begin{equation}\label{prop2eq1}
 0\leqslant\inf_{\psi\in J^{-1}\left(\sk{\mappingcone}\right)}\Tr\left(\Choimatr{\phi_{\lambda}}\Choimatr{\psi}\right)=\inf_{\psi\in J^{-1}\left(\sk{\mappingcone}\right)}\Tr\left(\left(\One-\lambda\Choimatr{\phi}\right)\Choimatr{\psi}\right)=1-\lambda\sup_{\psi\in J^{-1}\left(\sk{\mappingcone}\right)}\Tr\left(\Choimatr{\phi}\Choimatr{\psi}\right)=1-\lambda\norm{\phi}_{\mappingcone},
\end{equation}
iff $\norm{\phi}_{\mappingcone}\leqslant 1/\lambda$.
\end{proof}
\end{proposition}
In particular, if we apply the proposition to the cone of $k$-positive maps $\Kpositive{\kPmaps{k}}$ and use the fact that $\norm{\Choimatr{\phi}}_{S\left(k\right)}=\norm{\phi}_{\kPmaps{k}}$, we get the following result by Johnston and Kribs \cite{ref.JK09}.
\begin{corollary}\label{corr1}
 The map $\phi_{\lambda}=\Tr-\lambda\phi$, $\phi\in\CPmaps$, is $k$-positive iff
\begin{equation}\label{cor1eq1}
 \norm{\Choimatr{\phi}}_{S\left(k\right)}\leqslant\frac{1}{\lambda}.
\end{equation}
\end{corollary}
With the proper identification of norms, Proposition \ref{normphilambda} is a generalization of a characterization of Chruściński and Kossakowski \cite{ref.ChK09} of $k$-positivity in terms of Ky Fan norms. We first prove a lemma. Note that for a vector $\upsilon$, $\diad{\upsilon}$ denotes the rank one operator $\norm{\upsilon}^2\proj{\upsilon}$, where $\proj{\upsilon}$ is the one-dimensional projection onto $\setC\upsilon$.
\begin{lemma}\label{lemma1}
 Let $V=\sum_{i,j}V_{ij}e_{ij}\in\bh$, $\upsilon=\sum_{i,j}V_{ij}e_j\otimes e_i\in\hilbertspacetwo\otimes\hilbertspacetwo$, where $e_{ij}e_k=\delta_{jk}e_i$ for an orthonormal basis $e_1,\ldots,e_d$ for $\hilbertspacetwo$. Then we have
\begin{enumerate}[(i)]
 \item $\Choimatr{\Ad_V}=\diad{\upsilon}$
\item $\norm{\Choimatr{\Ad_V}}_{\textnormal{HS}}=\Tr\left(\Choimatr{\Ad_V}\right)=\norm{\upsilon
}^2=\norm{V}^2_{\textnormal{HS}}$
\item If $W\in\bh$ then
\begin{equation}\label{pointiii}
 \Tr\left(\Choimatr{\Ad_V}\Choimatr{\Ad_W}\right)=\abs{\innerpr{\upsilon}{\omega}}^2=\abs{\innerpr{V}{W}_{\textnormal{HS}}}^2,
\end{equation}
where $\omega$ is an analogue of $\upsilon$ for $W$. 
\end{enumerate}
\begin{proof}
 By the definition of $C_{\phi}$, we have
\begin{equation}\label{CAdV}
 C_{\Ad_V}=\sum_{i,j} e_{ij}\otimes Ve_{ij}\conj{V}=\sum_{i,j,m,n}V_{mi}\overline{V_{nj}} e_{ij}\otimes e_{mn}=\sum_{i,j,m,n}V_{mi}\overline{V_{nj}}\ddiad{e_i\otimes e_m}{e_j\otimes e_n}=\diad{\upsilon}.
\end{equation}
This proves $(i)$. To show $(ii)$, note that by $(i)$
\begin{equation}\label{CAdVsquare}
 \norm{\Choimatr{\Ad_V}}^2_{\textnormal{HS}}=\Tr\left(\Choimatr{\Ad_V}^2\right)=\norm{\upsilon}^2\Tr\left(\diad{\upsilon}\right)=\norm{\upsilon}^4=\Tr\left(\Choimatr{\Ad_V}\right)^2.
\end{equation}
Since furthermore $\norm{\upsilon}^2=\sum_{i,j}V_{ij}\overline{V_{ij}}=\norm{V}_{\textnormal{HS}}^2$, $(ii)$ follows. By $(i)$ 
\begin{equation}
 \Tr\left(\Choimatr{\Ad_V}\Choimatr{\Ad_W}\right)=\Tr\left(\diad{\upsilon}\diad{\omega}\right)=\abs{\innerpr{\upsilon}{\omega}}^2=\abs{\sum_{i,j}V_{ij}\overline{W_{ij}}}^2=\abs{\innerpr{V}{W}_{\textnormal{HS}}}^2,
\end{equation}
which proves $(iii)$.
\end{proof}
\end{lemma}
It turns out \cite{ref.ChK09,ref.JK09} that the Schmidt operator norms for maps of the form $\Ad_V$ can be effectively calculated in terms of the so-called {\it Ky Fan norms}.

\begin{definition}\label{KyFandef} Let $V$ be an element of $\bh$, $\dim\hilbertspacetwo=d$. For $k\in\left\{1,2,\ldots,d\right\}$, define the norm of $\kfnorm{k}{\Cdot}$ by
\begin{equation}\label{eqKyFandef}
 \norm{V}_{\left(k\right)}^2=\sum_{i=1}^k\sigma_i^2,
\end{equation}
where $\sigma_i^2$ is the $i$-th greatest eigenvalue of $V\conj{V}$. Then $\kfnorm{k}{\Cdot}^2$ is the $k$-th Ky Fan norm of $V\conj{V}$.
\end{definition}

Another way to define the norm $\kfnorm{k}{\Cdot}$ follows from the next proposition.

\begin{proposition}\label{KyFanchar}Take $V\in\bh$. Then
 \begin{equation}\label{KyFanalternativeeq}
  \norm{V}_{\left(k\right)}^2=\sup_{\rk F= k}\Tr\left(FV\conj{V}\right),
 \end{equation}
where $F$ runs over projections of dimension $k$.
\begin{proof}
 Can be found in \cite{ref.ChK09} or in \cite{ref.JK09}.
\end{proof}
\end{proposition}

\begin{theorem}\label{thm2notes}
 Let $V\in\bh$. Then 
\begin{equation}\label{KyFanofV}
 \kkfnorm{V}^2=\sup\left\{\Tr\left(\Choimatr{\Ad V}\Choimatr{\Ad W}\right)|\rk W\leqslant k,\Tr\left(\Choimatr{\Ad W}\right)=1\right\}=\sup_{\psi\in J^{-1}\left(\sk{\kPmaps{k}}\right)}\Tr\left(\Choimatr{{\Ad}_V}\Choimatr{\psi}\right)=\norm{{\Ad}_V}_{\kPmaps{k}}
\end{equation}
\begin{proof}
 The last equality in \eqref{KyFanofV} simply follows from the definition of the norm $\norm{.}_{\mappingcone}$, whereas the penultimate is a consequence of the fact that $J^{-1}\left(\sk{\kPmaps{k}}\right)=\convhull\left\{\Ad_W|\rk W\leqslant k,\Tr\left(\Choimatr{\Ad W}\right)=1\right\}$. In the remaining equality in \eqref{KyFanofV}, we first show ``$\leqslant$''. Let $W$ be as in \eqref{KyFanofV}. Since $W$ is of rank~$\leqslant k$, its range projection $E$ has dimension $\leqslant k$, and $W=EW$. We thus have by Lemma \ref{lemma1}, Proposition \ref{KyFanchar} and the Cauchy-Schwarz inequality
\begin{multline} 
 \Tr\left(\Choimatr{\Ad_V}\Choimatr{\Ad_W}\right)=\left|\hsprod{V}{W}\right|^2=\left|\Tr\left(V\conj{W}\right)\right|^2=\left|\Tr\left(V\conj{W}E\right)\right|^2\leqslant\Tr\left(EV\conj{\left(EV\right)}\right)\Tr\left(W\conj{W}\right)=\\=\Tr\left(EV\conj{V}\right)\hsnorm{W}^2\leqslant\sup_{\rk F\leqslant k}\Tr\left(FV\conj{V}\right)\cdot 1=\kfnorm{k}{V}^2.
\end{multline}
Now we prove~``$\geqslant$''. Since $\hilbertspacetwo$ is finite-dimensional, by compactness we can find a projection $E$ of dimension $\leqslant k$ such that $\kfnorm{k}{V}^2=\sup_{\rk F\leqslant k}\Tr\left(FV\conj{V}\right)=\Tr\left(EV\conj{V}\right)$. Take $W=EV/\kfnorm{k}{V}$. Then $\rk W\leqslant k$, and
\begin{equation}
 \hsnorm{W}^2=\frac{\Tr\left(EV\conj{\left(EV\right)}\right)}{\kfnorm{k}{V}^2}=\frac{\Tr\left(EV\conj{V}\right)}{\kfnorm{k}{V}^2}=1.
\end{equation}
In particular, $1=\Tr\left(W\conj{V}\right)/\kfnorm{k}{V}=\Tr\left(V\conj{W}\right)/\kfnorm{k}{V}$. By Lemma \ref{lemma1}, $(ii)$, $\Tr\left(\Choimatr{\Ad_W}\right)=1$. Furthermore,
\begin{multline}
 \kfnorm{k}{V}^2=\Tr\left(EV\conj{V}\right)=\kfnorm{k}{V}\Tr\left(V\conj{W}\right)=\kfnorm{k}{V}\Tr\left(V\conj{W}\right)\cdot 1=\kfnorm{k}{V}\Tr\left(V\conj{W}\right)\frac{\Tr\left(W\conj{V}\right)}{\kfnorm{k}{V}}=\\=\left|\Tr\left(V\conj{W}\right)\right|^2=\left|\hsprod{V}{W}\right|^2=\Tr\left(\Choimatr{\Ad_V}\Choimatr{\Ad_W}\right).
\end{multline}
Thus the $\sup$ is attained and we have the asserted equality.
\end{proof}
\end{theorem}

As an immediate corollary of Theorem \ref{thm2notes}, we get the result by Chru\'sci\'nski and Kossakowski \cite{ref.ChK09},
\begin{corollary}\label{corcharkpos}
 Let $V$ be an element of $\bh$. The map $\Tr-\lambda\Ad_V$ is $k$-positive if and only if
\begin{equation}\label{eqKyFankpos}
 \kfnorm{k}{V}^2\leqslant\frac{1}{\lambda}.
\end{equation}
\begin{proof}
 An immediate consequence of Proposition \ref{normphilambda} and Theorem \ref{thm2notes}.
\end{proof}
\end{corollary}
Int the special case $k=d=\dim\hilbertspacetwo$, we have $\hsnorm{V}^2=\kfnorm{d}{V}^2$. Since a map $\phi$ is $d$-positive iff $\phi$ is completely positive, we get
\begin{corollary}\label{corthree}
 Let $V\in\bh$. Then the map $\phi_{\lambda}=\Tr-\lambda\Ad_V$ is completely positive iff
\begin{equation}\label{Vhsle1overlambda}
 \hsnorm{V}^2\leqslant\frac{1}{\lambda}.
\end{equation}

\end{corollary}

It is instructive to derive the last two corollaries in yet another way. We first give a new proof of Corollary \ref{corthree}.
If $\hilbertspaceone=\hilbertspacetwo$, $\sk{\CPmaps}$ equals the set of states on $\bounded{\hilbertspacetwo\otimes\hilbertspacetwo}$. Indeed by \cite{ref.St09mappingcones}, $\CPmaps=\Kpositive{\CPmaps}$. Since a map $\phi$ belongs to $\CPmaps$ iff $\Choimatr{\phi}\geqslant 0$, and an operator $A$ is positive iff $\Tr\left(AB\right)\geqslant 0\,\forall_{B\geqslant 0}$, we have $\CPmaps=\dual{\CPmaps}=\dual{\Kpositive{\CPmaps}}$. Therefore $\sk{\CPmaps}$ is the set of states on $\bounded{\hilbertspacetwo\otimes\hilbertspacetwo}$ as asserted. Thus by Lemma \ref{lemma1},
\begin{equation}\label{kfnormAdhsnormV}
 \norm{{\Ad}_V}_{\CPmaps}=\sup_{\rho\textnormal{ - state}}\rho\left(\Choimatr{{\Ad}_V}\right)=\norm{C_{{\Ad}_V}}=\Tr\left(\Choimatr{{\Ad}_V}\right)=\hsnorm{V}^2,
\end{equation}
where the third equality follows since $\Choimatr{\Ad_V}$ is of rank $1$ and the fourth by $(ii)$ in Lemma \ref{lemma1}.

It turns out that the case of $k$-positive maps of the form $\Tr-\lambda\Ad_V$ can be reduced to the situation of Corollary \ref{corthree} using symmetries of the cone of $k$-positive maps. Let us start with the following lemma.

\begin{lemma}\label{lemma2Erling}
 Let $\mappingcone\in\Pmapsb{\hilbertspacetwo}$ be a symmetric mapping cone. Then a map $\phi\in\Pmapsb{\hilbertspacetwo}$ is $\mappingcone$-positive if and only if $\psi\circ\phi\in\CPmaps$ for all $\psi\in\dual{\mappingcone}$.
\begin{proof}
 We know from \cite[Thm. 2]{ref.St09mappingcones} that $\phi\in\Kpositive{\mappingcone}$ iff $\phi\in\mappingcone$. By \cite[Thm. 12]{ref.St09dual} and \cite[Cor. 9]{ref.St09mappingcones}, this is equivalent to saying that $\phi\in\dual{\Kpositive{\dual{\mappingcone}}}$. By \cite[Thm. 1]{ref.St09dual}, this is equivalent to the condition that $\psi\circ\phi\in\CPmaps$ for all $\psi\in\dual{\mappingcone}$.
\end{proof}
\end{lemma}

\begin{theorem}\label{theoremfour}
 Let $\phi\in\Pmapsb{\hilbertspacetwo}$, $k\in\left\{1,\ldots,d\right\}$. Then the following conditions are equivalent.
\begin{enumerate}[(i)]
\item $\phi$ is $k$-positive,
\item$\Ad_F\circ\,\phi\in\CPmaps$ for all  $k$-dimensional projections $F$ in $\hilbertspacetwo$,
\item$\phi\circ\Ad_F\in\CPmaps$ for all  $k$-dimensional projections $F$ in $\hilbertspacetwo$,
\item $\Ad_F\circ\,\phi\circ\Ad_E\in\CPmaps$ for all  $k$-dimensional projections $F$ and $E$ in $\hilbertspacetwo$,
\end{enumerate}
\begin{proof}
Since $\kPmaps{k}$ is a symmetric mapping cone, Lemma \ref{lemma2Erling} is applicable to $\kPmaps{k}$. Since $\kSPmapsb{k}{\hilbertspacetwo}$ is the dual of $\kPmaps{k}$ and is generated by maps $\Ad_F$ with $F$ projections of dimension $\leqslant k$, it follows that $(i)\Leftrightarrow(ii)$. Since $\kPmaps{k}$ and $\CPmaps$ are closed under taking adjoints, $(ii)\Leftrightarrow(iii)$. The equivalence $(i)\Leftrightarrow(iv)$ is a restatement of \cite[Lemma 2]{ref.St09tensorpowers}. 
\end{proof}
\end{theorem}
Similar characterization theorems have been discussed
in \cite{ref.SSZ09} and in a more general
 form in \cite{ref.St09mappingcones, ref.Sk10}. Using the theorem above, Proposition \ref{KyFanchar} and Corollary \ref{corthree}, one almost immediately gets Corollary \ref{corcharkpos}. 
Namely, for $\phi$ of the form $\phi_{\lambda}=\Tr-\lambda\Ad_V$ and $E$, $F$ projections of dimension $k$, we have the map $\phi^{EF}_{\lambda}$,
\begin{equation}\label{phiEFdef}
 \phi^{EF}_{\lambda}\left(x\right):={\Ad}_F\circ\phi_{\lambda}\circ{\Ad}_E\left(x\right)=F\Tr\left(ExE\right)-\lambda{\Ad}_{FV}\left(ExE\right).
\end{equation}
Let $U$ be a unitary operator such that $E=UF\conj{U}$. One can rewrite \eqref{phiEFdef} as
\begin{equation}\label{phiEFsecondform}
 \phi_{\lambda}^{EF}\left(x\right)=F\Tr\left(UF\conj{U}xUF\conj{U}\right)-\lambda{\Ad}_{FV}\left(UF\conj{U}xUF\conj{U}\right)=F\Tr\left(F\conj{U}xUF\right)-\lambda{\Ad}_{FVU}\left(F\conj{U}xUF\right).
\end{equation}
Since the map $x\mapsto\conj{U}xU$ is an isomorphism, $\phi^{EF}$ is completely positive iff the map $\phi'_{\lambda}$ defined by
\begin{equation}\label{phiprimedef}
\phi'_{\lambda}\left(y\right)=F\Tr\left(FyF\right)-\lambda{\Ad}_{FVU}\left(FyF\right)
\end{equation}
is a completely positive map of $\bounded{F\hilbertspacetwo}$ into itself. By Corollary \ref{corthree}, this happens iff
\begin{equation}\label{FVU}
\hsnorm{FVU}^2\leqslant\frac{1}{\lambda} 
\end{equation}
or equivalently
\begin{equation}\label{FVUorequation}
 \frac{1}{\lambda}\geqslant\Tr\left(FVU\conj{U}\conj{V}F\right)=\Tr\left(FV\conj{V}F\right)=\Tr\left(FV\conj{V}\right).
\end{equation}
One thing which needs a comment in equations \eqref{FVU} and \eqref{FVUorequation} is that the trace of an operator in $\bounded{F\hilbertspacetwo}$ equals $\Tr\left(F\Cdot F\right)$, where $\Tr$ denotes the usual trace in $\hilbertspacetwo$. To see that this is the case, it is enough to calculate $\Tr$ in a basis of $\hilbertspacetwo$ where the first $k$ vectors belong to $F\hilbertspacetwo$ and the remaining $d-k$ ones to $F\hilbertspacetwo^{\bot}$.
Taking a supremum over $k$-dimensional projections $F$ in \eqref{FVUorequation}, by Theorem~\ref{theoremfour}, we get
\begin{equation}\label{KyFanfinally}
 \phi_{\lambda}\in\kPmaps{k}\,\Leftrightarrow\,\kfnorm{k}{V}^2=\sup_{\rk F=k}\Tr\left(FV\conj{V}\right)\leqslant\frac{1}{\lambda},
\end{equation}
where we also used Proposition \ref{KyFanchar}. Formula \eqref{KyFanfinally} is the same as in Corollary \ref{corcharkpos}. Again, this is the result by Chru\'sci\'nski and Kossakowski on $k$-positive maps of the form $\Tr-\lambda\Ad_V$ \cite{ref.ChK09}.

\section{Completely and $\mappingcone$-entangled subspaces}
It turns out that if $\phi$ is a positive map then the negative part $\Choimatrminus{\phi}$ of the Choi matrix $\Choimatr{\phi}$ of $\phi$ contains much interesting information. We shall study this observation in the present section.

\begin{definition}\label{def1}
Let $\mappingcone$ be a mapping cone on $\hilbertspacetwo$, $\mappingcone\not\subset\CPmapsb{\hilbertspacetwo}$ and let $\hilbertspaceone$ be another Hilbert space. As before, 
\begin{equation}\label{skskdef}
 \sk{\mappingcone}=\left\{\rho\in\boundedconj{\hilbertspaceone\otimes\hilbertspacetwo}\vline\rho=\Tr\left(\Choimatr{\psi}.\right)\textnormal{ is a state },\psi\in\dual{\Kpositive{\mappingcone}}\right\}
\end{equation}
and $\Kpositive{\mappingcone}$ denotes the set of $\mappingcone$-positive maps $\bk\rightarrow\bh$. We say a state $\omega$ on $\bounded{\hilbertspaceone\otimes\hilbertspacetwo}$ is {\it $\mappingcone$-entangled} if $\omega\not\in\sk{\mappingcone}$.
\end{definition}
Note that if $\mappingcone=\Pmapsb{\hilbertspacetwo}$ then $\omega$ is $\mappingcone$-entangled iff $\omega$ is entangled since $\sk{\Pmapsb{\hilbertspacetwo}}$ is the set of separable states. We can now state the main result of this section. Note that the projection $e$ in the theorem will be the support projection for $\Choimatrminus{\phi}$, or equivalently, the range projection of $\Choimatrminus{\phi}$.
\begin{theorem}\label{mainthmsec2}
 Let $e$ be a projection in $\bounded{\hilbertspaceone\otimes\hilbertspacetwo}$ and $\mappingcone$ a mapping cone on $\hilbertspacetwo$ with $\mappingcone\not\subset\CPmapsb{\hilbertspacetwo}$. Then each state $\omega$ on $\bounded{\hilbertspaceone\otimes\hilbertspacetwo}$ with support in $e$ is $\mappingcone$-entangled iff there exists a $\mappingcone$-positive map $\phi:\bk\rightarrow\bh$ with support $\Choimatrminus{\phi}=e$.
\begin{proof}
 Suppose $\supp\Choimatrminus{\phi}=e$. If $\omega$ is a state with $\supp\omega\leqslant e$ then
\begin{equation}\label{eqomegaCphi}
 \omega\left(\Choimatr{\phi}\right)=\omega\left(e\Choimatr{\phi}e\right)=-\omega\left(\Choimatrminus{\phi}\right)<0,
\end{equation}
Thus if $\omega=\Tr\left(\Choimatr{\psi}\Cdot\right)$, then $\Choimatr{\psi}\not\in\dual{\Kpositive{\mappingcone}}$, hence $\omega$ is $\mappingcone$-entangled. To show the converse, let $\mu=\sup_{\rho\in\sk{\mappingcone}}\rho\left(e\right)$. We claim that $\mu<1$. Indeed, $1=\norm{e}=\sup\left\{\Tr\left(eh\right)\vline 0\leqslant h\leqslant 1,\Tr\left(h\right)=1\right\}$. Now $\Tr\left(eh\right)=\Tr\left(h\right)=1$ iff $h\leqslant e$. Thus the state $\Tr\left(h\Cdot\right)$ is by assumption $\mappingcone$-entangled, hence $\Tr\left(h\Cdot\right)\not\in\sk{\mappingcone}$. Thus $\Tr\left(e\Choimatr{\psi}\right)<1$ for all states $\Tr\left(\Choimatr{\psi}.\right)\in\sk{\mappingcone}$. By compactness of $\sk{\mappingcone}$ and continuity of the maps $\psi\mapsto\Tr\left(e\Choimatr{\psi}\right)$, $\mu<1$, as claimed. Let $\lambda=1/\mu$, and $\phi$ be defined by $\Choimatr{\phi}=\One-\lambda e$. If $\Tr\left(\Choimatr{\psi}\Cdot\right)\in\sk{\mappingcone}$, then 
\begin{equation}
 \Tr\left(\Choimatr{\phi}\Choimatr{\psi}\right)=\One-\lambda\Tr\left(e\Choimatr{\psi}\right)\geqslant 1-\lambda\mu=0.
\end{equation}
Thus $\phi\in\ddual{\Kpositive{\mappingcone}}$. By Theorem 6 in \cite{ref.St09dual}, $\ddual{\Kpositive{\mappingcone}}=\Kpositive{\mappingcone}$, so $\phi$ is $\mappingcone$-positive. Since $\Choimatr{\phi}=\left(\One-e\right)-\left(\lambda-1\right)e$, $\Choimatrplus{\phi}=\One-e$, $\Choimatrminus{\phi}=\left(\lambda-1\right)e$, which has support $e$.
\end{proof}
\end{theorem}
As an immediate corollary we have
\begin{corollary}
 If $\phi:\bk\rightarrow\bh$ is $\mappingcone$-positive with $\mappingcone$ as in Theorem \ref{mainthmsec2}, then every state with support in the support of $\Choimatrminus{\phi}$ is $\mappingcone$-entangled.
\end{corollary}

Following \cite{ref.Partha04} we say that a projection $e$ in $\bounded{\hilbertspaceone\otimes\hilbertspacetwo}$ is {\it completely entangled} if each state $\omega$ with support in $e$ is entangled. In the special case when $\mappingcone=\Pmapsb{\hilbertspaceone}$, we get
\begin{corollary}\label{cor2sec2}
 Let $e$ be a projection in $\bounded{\hilbertspaceone\otimes\hilbertspacetwo}$. Then $e$ is completely entangled iff there exists a positive map $\phi:\bk\rightarrow\bh$ such that $\supp\Choimatrminus{\phi}=e$.
\end{corollary}
It is natural to ask what limits on the dimension of the support of $\Choimatrminus{\phi}$ are implied by the fact that $\phi$ belongs to a mapping cone $\mappingcone$. For the cones of $k$-positive maps, the question has received some attention in the past \cite{ref.Timoney00,ref.J04}. Another estimate follows from Theorem \ref{mainthmsec2} and the results of \cite{ref.CMW08}. 
\begin{corollary}\label{cor3sec2}
 Let $\phi:\bk\rightarrow\bh$ be $k$-positive, $k\in\left\{1,2,\ldots,d-1\right\}$, $d\leqslant\min\left(\dim\hilbertspaceone,\hilbertspacetwo\right)$. Then $\dim\supp\Choimatrminus{\phi}\leqslant\left(m-k\right)\left(n-k\right)$, where $m=\dim\hilbertspaceone$, $n=\dim\hilbertspacetwo$.
\begin{proof}
 Since $\dual{\kPmaps{k}}=\kSPmaps{k}$ the vectors in the support of $\Choimatrminus{\phi}$
 have Schmidt rank $\geqslant k+1$ (since $\sk{\kPmaps{k}}$ consists of states with density operators in $\kSPmaps{k}$). Then by a result of \cite[Thm. 11]{ref.CMW08}, this means that
\begin{equation}
 \dim\supp\Choimatrminus{\phi}\leqslant\left(m-\left(k+1\right)+1\right)\left(n-\left(k+1\right)+1\right)=\left(m-k\right)\left(n-k\right).
\end{equation}
\end{proof}
\end{corollary}
Note that the same estimate was earlier obtained in \cite[Prop. 2]{ref.Gniewko08}.

\section{Tensor products}
We now turn to an example about $k$-positivity of tensor products of positive maps. Throughout the section we will be using the elementary fact that $\Choimatr{\phi^{\otimes n}}=\Choimatr{\phi}^{\otimes n}$. We know from Corollary \ref{corr1} and the equality $\Choimatr{\Tr}=\One$ that estimation of the operator Schmidt norm $\left\|A\right\|_{S\left(k\right)}$ for $k\in\mathbbm{N}$ is crucial for checking whether a Choi matrix of the form $\One-\lambda A$ corresponds by the Jamiołkowski-Choi isomorphism to a $k$-positive map or not. Explicitly, the norm $\norm{A}_{S\left(k\right)}$ for an operator $A\in\bh$ is defined \cite{ref.JK09} by the formula
\begin{equation}\label{Schmidtnorm}
 \left\|A\right\|_{S\left(k\right)}=\sup_{\textrm{SR}\left(\psi\right)\leqslant k}\left<\psi|A|\psi\right>,
\end{equation}
where the supremum is taken over vectors $\psi$ of Schmidt rank $\leqslant k$, i.e. vectors $\psi=\sum_{i,j}\psi_{ij}e_i\otimes e_j$ in $\hilbertspacetwo\otimes\hilbertspacetwo$ such that the matrix $\left[\psi_{ij}\right]$ is of rank $\leqslant k$.

 When distillability of quantum states is discussed, the case $k=2$ is of particular interest. More precisely, the question of distillability of a bipartite state $\rho$ is equivalent to the question whether for some $n\in\mathbbm{N}$ the matrix $\left(\left(\One\otimes\transpose\right)\rho\right)^{\otimes n}$ corresponds by the Choi-Jamiołkowski isomorphism to a $2$-positive map or not. When $\left(\left(\One\otimes\transpose\right)\rho\right)^{\otimes n}$ corresponds to a $2$-positive map for some $n\in\mathbbm{N}$, we say that $\rho$ is $n$-copy undistillable \cite{ref.DiVicenzo00,ref.DCLB00}. Of special importance are density matrices such that $\left(\One\otimes\transpose\right)\rho$ is proportional to $\One-\lambda p$, with $p$ a projection, especially when $p$ projects onto $\psi_+=\frac{1}{\sqrt{d}}\sum_{i=1}^de_i\otimes e_i$, the maximally entangled state. We know from the previous discussion (cf. Lemma \ref{lemma1} plus the equality  $\Choimatr{\Tr}=\One$ in Proposition \ref{normphilambda}) that operators of the form $\One-\lambda p$ correspond to maps of the form $\Tr-\lambda\Ad_V$. Thus, the question of distillability is about $2$-positivity of maps $\left(\Tr-\lambda \Ad_V\right)^{\otimes n}$. For $n=1$, the question has fully been solved in \cite{ref.ChK09} and the answer can be found in Corollary \ref{corcharkpos} above. For $n=2$, we have
\begin{equation}
 \Choimatr{\left(\Tr-\lambda \Ad_V\right)^{\otimes 2}}=\left(\One-\lambda p\right)^{\otimes 2}=\One\otimes\One-\lambda\left(\One\otimes p+p\otimes\One\right)+\lambda^2 p\otimes p.
\end{equation}
Thus
\begin{equation}
 \Choimatr{\left(\Tr-\lambda \Ad_V\right)^{\otimes 2}}\geqslant\One\otimes\One-\lambda\left(\One\otimes p+p\otimes\One\right)
\end{equation}
and any sufficient $2$-positivity criterion for the map $J^{-1}\left(\One\otimes\One-\lambda\left(\One\otimes p+p\otimes\One\right)\right)$ will work for $\left(\Tr-\lambda \Ad_V\right)^{\otimes 2}$ as well. But $\One\otimes\One-\lambda\left(\One\otimes p+p\otimes\One\right)$ is of the form $\One-\lambda A$, so by \cite{ref.JK09}, or Corollary \ref{corr1} above, we have
\begin{proposition}\label{prop1}
 If $\frac{1}{\lambda}\geqslant\left\|\One\otimes p+p\otimes\One\right\|_{S\left(2\right)}$, then the map $\left(\Tr-\lambda \Ad_V\right)^{\otimes 2}$ is $2$-positive.
\end{proposition}
Obviously, we have
$\left\|\One\otimes p+p\otimes\One\right\|_{S\left(2\right)}\leqslant 2\left\|\One\otimes p\right\|_{S\left(2\right)}$.
Let $p$ be a projection onto a vector $\psi$ of norm one. To make further estimates, we need to introduce the concept of singular values and Schmidt vector norms \cite{ref.JK09} for $\psi$.
\begin{definition}\label{defsingularvaluespsi}
By the $i$-th singular value $\sigma_i\left(\psi\right)$ for a vector $\psi=\sum_{i,j}\psi_{ij}e_i\otimes e_j\in\hilbertspacetwo\otimes\hilbertspacetwo$ we mean the $i$-th greatest singular value of the coordinate matrix $\left[\psi_{ij}\right]$.
\end{definition}
\begin{definition}\label{defSchmidtvectornorm}
 The $k$-th Schmidt vector norm $\norm{\Cdot}_{s\left(k\right)}$ for a vector $\psi=\sum_{i,j}\psi_{ij}e_i\otimes e_j\in\hilbertspacetwo\otimes\hilbertspacetwo$ is the norm $\kfnorm{k}{\left[\psi_{ij}\right]}$ of the coordinate matrix $\left[\psi_{ij}\right]$ (cf. Definition \ref{KyFandef}).
\end{definition}
The eigenvectors of $\One\otimes p$ corresponding to the eigenvalue $1$ are of the form $v_i\otimes\psi$, where $\left\{v_i\right\}_{j=1}^d$ can be any orthonormal basis of $\hilbertspacetwo$. Other eigenvectors correspond to the eigenvalue $0$ and we can neglect them in our discussion.
We have the following lemma concerning the singular values of a tensor product of two vectors.
\begin{lemma}\label{lemmasingproduct}
 Let $\phi\otimes\psi$ be a vector in $\left(\hilbertspacetwo\hat\otimes\hilbertspacetwo\right)\otimes\left(\hilbertspacetwo\hat\otimes\hilbertspacetwo\right)$, where $\phi,\psi\in\hilbertspacetwo\hat\otimes\hilbertspacetwo$ and $\hat\,$ is used to distinguish between the two tensor products that we are using. The singular values of $\phi\otimes\psi$ with respect to the tensor product $\hat\otimes$ are equal to $\sigma_i\left(\phi\right)\sigma_j\left(\psi\right)$.
\begin{proof}
 Let $\phi=\sum_{i,j}\phi_{ij}e_i\hat\otimes e_j$, $\psi=\sum_{k,l}\psi_{kl}e_k\hat\otimes e_l$ be the decompositions of $\phi$, $\psi$ (resp.) in an orthonormal basis of $\hilbertspacetwo\hat\otimes\hilbertspacetwo$. The decomposition of $\phi\otimes\psi$ with respect to an orthonormal basis of $\left(\hilbertspacetwo\hat\otimes\hilbertspacetwo\right)\otimes\left(\hilbertspacetwo\hat\otimes\hilbertspacetwo\right)$ is
\begin{equation}
 \sum_{i,j,k,l}\phi_{ij}\psi_{kl}\left(e_i\otimes e_k\right)\hat\otimes\left(e_j\otimes e_l\right).
\end{equation}
Thus the singular values for $\phi\otimes\psi$ are the singular values of the matrix $\mathbbm{A}=\left[A_{ik,jl}\right]$, where $A_{ik,jl}=\phi_{ij}\psi_{kl}$. But $\mathbbm{A}=\Phi\otimes\Psi$, where $\Phi=\left[\phi_{ij}\right]$, $\Psi=\left[\psi_{kl}\right]$. Thus $\mathbbm{A\conj{A}}=\Phi\conj{\Phi}\otimes\Psi\conj{\Psi}$ and the eigenvalues for $\mathbbm{A\conj{A}}$ are products of $\sigma^2_i\left(\phi\right)$ and $\sigma^2_j\left(\psi\right)$ ($i,j\in\left\{1,2,\ldots,d\right\}$). In other words, the singular values for $\phi\otimes\psi$ are equal to the products of the singular values for $\phi$ and for $\psi$.
\end{proof}
\end{lemma}
In short, the singular values of the vectors $v_i\otimes\psi$ are of the form $\sigma_k\left(v_i\right)\sigma_l\left(\psi\right)$, $l,m=1,2,\ldots,d$.  Thus we have
\begin{equation}\label{est2}
 \sigma_2\left(v_i\otimes\psi\right)\leqslant\sigma_1\left(v_i\otimes\psi\right)\leqslant\sigma_1\left(v_i\right)\sigma_1\left(\psi\right)
\end{equation}
and
\begin{equation}\label{est3}
\norm{v_i\otimes\psi}_{s\left(2\right)}^2\leqslant 2\sigma_1\left(v_i\right)^2\sigma_1\left(\psi\right)^2=2\norm{v_i}_{s\left(1\right)}^2\norm{\psi}_{s\left(1\right)}^2
\end{equation}
By \cite{ref.JK09}, Proposition 4.8, we have
\begin{equation}\label{est4}
 \left\|\One\otimes p\right\|_{S\left(2\right)}\leqslant\sum_{i=1}^d\left\|v_i\otimes\psi\right\|_{s\left(2\right)}^2\leqslant 2\left\|\psi\right\|_{s\left(1\right)}^2\sum_{i=1}^d\left\|v_i\right\|_{s\left(1\right)}^2
\end{equation}
The vectors $v_i$ can be chosen as an arbitrary orthonormal basis of the space in question (of dimensionality $d$). In particular, from \cite{ref.Werner01} we know that for arbitrary $d$, it is possible to construct an orthonormal basis of maximally entangled states. In such case $\left\|v_i\right\|_{s\left(1\right)}^2=\frac{1}{d}$ for all $i$, which is the minimum value that $\left\|\Cdot\right\|_{s\left(1\right)}$ can take in general, so optimizing the choice of $v_i$ in \eqref{est4} gives
\begin{equation}\label{est5}
 \left\|\One\otimes p\right\|_{S\left(k\right)}\leqslant 2\left\|\psi\right\|_{s\left(1\right)}^2\cdot d\cdot\frac{1}{d}=2\left\|\psi\right\|_{s\left(1\right)}^2
\end{equation}
Now we only need to plug this into $\left\|\One\otimes p+p\otimes\One\right\|_{S\left(2\right)}\leqslant 2\left\|\One\otimes p\right\|_{S\left(2\right)}$ to get
\begin{equation}
 \left\|\One\otimes p+p\otimes\One\right\|_{S\left(k\right)}\leqslant 4\left\|\psi\right\|_{s\left(1\right)}^2=4\left\|p\right\|_{S\left(1\right)}
\end{equation}
We can use this result to give a concrete estimate in Proposition \ref{prop1}.
\begin{proposition}\label{prop2}
 If $\frac{1}{\lambda}\geqslant 4\left\|p\right\|_{S\left(1\right)}$, $\Choimatr{\Ad_V}=p$, the map $\left(\Tr-\lambda \Ad_V\right)^{\otimes 2}$ is $2$-positive.
\end{proposition}
Furthermore, by Theorem \ref{mainthmsec2} we have that the support of $\One\otimes p+p\otimes\One$ is $\kPmaps{2}$-entangled. We also make the following observation,
\begin{corollary}\label{cor1}
 If $\Tr-4\lambda\Ad_V$ with $\hsnorm{V}=1$ is $1$-positive, then $\left(\Tr-\lambda \Ad_V\right)^{\otimes 2}$ is $2$-positive.
\begin{proof}
 By Corollary \ref{corr1}, we know that $\phi_{4\lambda}=\Tr-4\lambda\Ad_V$ is $1$-positive iff $\norm{p}_{S\left(1\right)}\leqslant \frac{1}{4\lambda}$ iff $\frac{1}{\lambda}\geqslant 4\norm{p}_{S\left(1\right)}$. This is precisely the condition in Proposition \ref{prop2}.
\end{proof}
\end{corollary}
Since $\left\|p\right\|_{S\left(1\right)}\geqslant\frac{1}{d}$ and we are only interested in $\lambda>1$ (otherwise, $\Tr-\lambda\Ad_V\in\CPmaps$), Proposition \ref{prop2} can only be useful when $d>4$.


\section{Acknowledgement}
The authors would like to thank Dariusz Chruściński, Andrzej Kossakowski, Marcin Marciniak and Karol \.Zyczkowski for comments on the manuscript. The project was operated within the Foundation
for Polish Science International Ph.D. Projects Programme co-financed
by the European Regional Development Fund covering, under the agreement
no. MPD/2009/6, the Jagiellonian University International Ph.D. Studies in
Physics of Complex Systems. Some of the results presented in the paper came to existence during a visit of one of the authors (Ł.S.) to the University of Oslo, where he enjoyed hospitality at the Mathematics Institute and was financially supported by Scholarschip and Training Fund, operated by Foundation for the Development of the Education System.

\bibliographystyle{elsarticle-num-names}
\bibliography{paper_kpositivity_JFA}

\begin{thebibliography}{21}
\providecommand{\natexlab}[1]{#1}
\providecommand{\url}[1]{\texttt{#1}}
\providecommand{\urlprefix}{URL }
\expandafter\ifx\csname urlstyle\endcsname\relax
  \providecommand{\doi}[1]{doi:\discretionary{}{}{}#1}\else
  \providecommand{\doi}[1]{doi:\discretionary{}{}{}\begingroup
  \urlstyle{rm}\url{#1}\endgroup}\fi
\providecommand{\bibinfo}[2]{#2}

\bibitem[{Chru{\'s}ci{\'n}ski and Kossakowski(2009)}]{ref.ChK09}
\bibinfo{author}{D.~Chru{\'s}ci{\'n}ski}, \bibinfo{author}{A.~Kossakowski},
  \bibinfo{title}{Spectral conditions for positive maps},
  \bibinfo{journal}{Commun. Math. Phys.} \bibinfo{volume}{290}
  (\bibinfo{year}{2009}) \bibinfo{pages}{1051--1064}.

\bibitem[{DiVincenzo et~al.(2000)DiVincenzo, Shor, Smolin, Terhal, and
  Thapliyal}]{ref.DiVicenzo00}
\bibinfo{author}{D.~DiVincenzo}, \bibinfo{author}{P.~Shor},
  \bibinfo{author}{J.~Smolin}, \bibinfo{author}{B.~Terhal},
  \bibinfo{author}{A.~Thapliyal}, \bibinfo{title}{Evidence for bound entangled
  states with negative partial transpose}, \bibinfo{journal}{Phys. Rev. A}
  \bibinfo{volume}{61} (\bibinfo{year}{2000}) \bibinfo{pages}{062312}.

\bibitem[{Choi(1975)}]{ref.Choi75}
\bibinfo{author}{M.-D. Choi}, \bibinfo{title}{Completely positive linear maps
  on complex matrices}, \bibinfo{journal}{Lin. Alg. Appl.} \bibinfo{volume}{10}
  (\bibinfo{year}{1975}) \bibinfo{pages}{285--290}.

\bibitem[{Kraus(1971)}]{ref.Kraus}
\bibinfo{author}{K.~Kraus}, \bibinfo{title}{General state changes in quantum
  theory}, \bibinfo{journal}{Ann. Phys.} \bibinfo{volume}{64}
  (\bibinfo{year}{1971}) \bibinfo{pages}{311--335}.

\bibitem[{Ando(2004)}]{ref.Ando04}
\bibinfo{author}{T.~Ando}, \bibinfo{title}{Cones and norms in the tensor
  product of matrix spaces}, \bibinfo{journal}{Lin. Alg. Appl.}
  \bibinfo{volume}{379} (\bibinfo{year}{2004}) \bibinfo{pages}{3--41}.

\bibitem[{Horodecki et~al.(2003)Horodecki, Shor, and Ruskai}]{ref.HSR03}
\bibinfo{author}{M.~Horodecki}, \bibinfo{author}{P.~Shor},
  \bibinfo{author}{M.~Ruskai}, \bibinfo{title}{Entanglement breaking channels},
  \bibinfo{journal}{Rev. Math. Phys.} \bibinfo{volume}{15}
  (\bibinfo{year}{2003}) \bibinfo{pages}{629--641}.

\bibitem[{St{\o}rmer(1986)}]{ref.St86}
\bibinfo{author}{E.~St{\o}rmer}, \bibinfo{title}{Extension of positive maps
  into $\bh$}, \bibinfo{journal}{J. Funct. Anal.} \bibinfo{volume}{66}
  (\bibinfo{year}{1986}) \bibinfo{pages}{235--254}.

\bibitem[{St{\o}rmer(????)}]{ref.St09mappingcones}
\bibinfo{author}{E.~St{\o}rmer}, \bibinfo{title}{Mapping cones of positive
  maps}, \bibinfo{journal}{Math. Scand.} \bibinfo{note}{To appear, preprint
  arXiv:0906.0472}.

\bibitem[{Jamio{\l}kowski(1972)}]{ref.J72}
\bibinfo{author}{A.~Jamio{\l}kowski}, \bibinfo{title}{Linear transformations
  which preserve trace and positive semidefinitness of operators},
  \bibinfo{journal}{Rep. Math. Phys.} \bibinfo{volume}{3}
  (\bibinfo{year}{1972}) \bibinfo{pages}{275--278}.

\bibitem[{Skowronek et~al.(2009)Skowronek, St{\o}rmer, and
  {\.Z}yczkowski}]{ref.SSZ09}
\bibinfo{author}{{\L}.~Skowronek}, \bibinfo{author}{E.~St{\o}rmer},
  \bibinfo{author}{K.~{\.Z}yczkowski}, \bibinfo{title}{Cones of positive maps
  and their duality relations}, \bibinfo{journal}{J. Math. Phys.}
  \bibinfo{volume}{50} (\bibinfo{year}{2009}) \bibinfo{pages}{062106}.

\bibitem[{St{\o}rmer(2009{\natexlab{a}})}]{ref.St09dual}
\bibinfo{author}{E.~St{\o}rmer}, \bibinfo{title}{Duality of cones of positive
  maps}, \bibinfo{journal}{M{\"u}nster J. Math.} \bibinfo{volume}{2}
  (\bibinfo{year}{2009}{\natexlab{a}}) \bibinfo{pages}{299--310}.

\bibitem[{Johnston and Kribs(2010)}]{ref.JK09}
\bibinfo{author}{N.~Johnston}, \bibinfo{author}{D.~Kribs}, \bibinfo{title}{A
  {F}amily of {N}orms {W}ith {A}pplications {I}n {Q}uantum {I}nformation
  {T}heory}, \bibinfo{journal}{J. Math. Phys.} \bibinfo{note}{To appear,
  preprint arXiv:0909.3907v3}.

\bibitem[{St{\o}rmer(2009{\natexlab{b}})}]{ref.St09tensorpowers}
\bibinfo{author}{E.~St{\o}rmer}, \bibinfo{title}{Tensor powers of $2$-positive
  maps}, \bibinfo{note}{to appear}, \bibinfo{year}{2009}{\natexlab{b}}.

\bibitem[{Skowronek(2010)}]{ref.Sk10}
\bibinfo{author}{{\L}.~Skowronek}, \bibinfo{title}{Theory of {G}eneralized
  {M}apping {C}ones in the {F}inite-{D}imensional {C}ase},
  \bibinfo{note}{preprint}, \bibinfo{year}{2010}.

\bibitem[{Parthasarathy(2004)}]{ref.Partha04}
\bibinfo{author}{K.~Parthasarathy}, \bibinfo{title}{On the maximal dimension of
  a completely entangled subspace for finite level quantum systems},
  \bibinfo{journal}{Proc. Indian Acad. Sci. (Math. Sci.)} \bibinfo{volume}{114}
  (\bibinfo{year}{2004}) \bibinfo{pages}{365--374}.

\bibitem[{Timoney(2000)}]{ref.Timoney00}
\bibinfo{author}{R.~Timoney}, \bibinfo{title}{A note on positivity of
  elementary operators}, \bibinfo{journal}{Bull. London Math. Soc.}
  \bibinfo{volume}{32} (\bibinfo{year}{2000}) \bibinfo{pages}{229--234}.

\bibitem[{Jamio{\l}kowski(2004)}]{ref.J04}
\bibinfo{author}{A.~Jamio{\l}kowski}, \bibinfo{title}{Some Remarks on the Role
  of Minimal Length of Positive Maps in Constructing Entanglement Witnesses},
  \bibinfo{journal}{Open Systems \& Information Dynamics} \bibinfo{volume}{11}
  (\bibinfo{year}{2004}) \bibinfo{pages}{385--390}.

\bibitem[{Cubitt et~al.(2008)Cubitt, Montanaro, and Winter}]{ref.CMW08}
\bibinfo{author}{T.~Cubitt}, \bibinfo{author}{A.~Montanaro},
  \bibinfo{author}{A.~Winter}, \bibinfo{title}{On the dimension of subspaces
  with bounded {S}chmidt rank}, \bibinfo{journal}{J. Math. Phys.}
  \bibinfo{volume}{49} (\bibinfo{year}{2008}) \bibinfo{pages}{022107}.

\bibitem[{Sarbicki(2008)}]{ref.Gniewko08}
\bibinfo{author}{G.~Sarbicki}, \bibinfo{title}{Spectral properties of
  entanglement witnesses}, \bibinfo{journal}{J. Phys. A} \bibinfo{volume}{41}
  (\bibinfo{year}{2008}) \bibinfo{pages}{375303}.

\bibitem[{D{\"u}r et~al.(2000)D{\"u}r, Cirac, Lewenstein, and
  Bruss}]{ref.DCLB00}
\bibinfo{author}{W.~D{\"u}r}, \bibinfo{author}{I.~Cirac},
  \bibinfo{author}{M.~Lewenstein}, \bibinfo{author}{D.~Bruss},
  \bibinfo{title}{Distillability and partial transposition in bipartite
  systems}, \bibinfo{journal}{Phys. Rev. A} \bibinfo{volume}{61}
  (\bibinfo{year}{2000}) \bibinfo{pages}{062313}.

\bibitem[{Werner(2001)}]{ref.Werner01}
\bibinfo{author}{R.~Werner}, \bibinfo{title}{All teleportation and dense coding
  schemes}, \bibinfo{journal}{J. Phys. A} \bibinfo{volume}{34}
  (\bibinfo{year}{2001}) \bibinfo{pages}{7081}.

\end{thebibliography}

\end{document}